\documentclass{amsproc}
\usepackage{graphicx}
\usepackage{amscd}
\usepackage{amsmath}
\usepackage{amsfonts}
\usepackage{amssymb}
\theoremstyle{plain}

\newtheorem{definition}{Definition}

\newtheorem{problem}{Problem}

\newtheorem{remark}{Remark}

\newtheorem{theorem}{Theorem}
\numberwithin{equation}{section}

\begin{document}
\title[Hereditarily Approximation Property]{An example of a Banach Space with a Subsymmetric Basis, which has the
Hereditarily Approximation Property}
\author{Eugene Tokarev}
\address{B.E. Ukrecolan, 33-81 Iskrinskaya str., 61005, Kharkiv-5, Ukraine}
\email{tokarev@univer.kharkov.ua}
\thanks{This paper is in final form and no version of it will be submitted for
publication elsewhere.}
\dedicatory{Dedicated to the memory of S. Banach.}
\begin{abstract}W.B. Johnson has constructed a series of Banach spaces non isomorphic to the
Hilbert one that have the hereditarily approximation property (shortly
hereditarily AP): all their subspaces also have the AP. All these examples
were ''sufficiently'' non symmetric and this fact allows Johnson to ask:
whether there exists any Banach space $X$ with symmetric (or, at least,
subsymmetric) basis, distinct from the Hilbert space such that each its
subspace has the AP? In this paper is shown that there is a Banach space with
a subsymmetric basis (non equivalent to any symmetric one), which enjoys the
hereditarily AP.
\end{abstract}
\maketitle

\section{Introduction}

A Banach space $X$ is said to have the \textit{approximation property}%
\textbf{\ (}shortly\textbf{: }AP) provided for every compact subset $K\subset
X$ and every $\varepsilon>0$ there exists a (bounded linear) operator
$u:X\rightarrow X$ with a finite dimensional range such that $\left\|
ux-x\right\|  \leq\varepsilon$ for all $x\in K$.

It will be said that $X$ has the \textit{hereditarily approximation property}
(\textit{hereditarily AP}) if every subspace of $X$ has the AP too. Certainly,
any Hilbert space has this property.

The long-standing approximation problem was\textit{: }

\textit{Whether every Banach space has the }AP?

It was solved in negative by Enflo [1], who constructed the first example of a
Banach space $X$, which$\ $does not have the AP. This result was strengthened
by A.~Szankowski [2], who showed that if there exists a such real $p\neq2$
that $l_{p}$ is finitely representable in $X$ (definition of this and other
notions will be given below) then $X$ contains a subspace $Y$, which does not
enjoy the AP. O.V. Reinov [3] constructed an example of a Banach space, that
has not the AP, with the property: $l_{p}$ is finitely representable in $X$
only for $p=2$.

Nevertheless, the natural hypothesis: \textit{''a Banach space }%
$X$\textit{\ has the hereditarily AP if and only if }$X$\textit{\ is
isomorphic to the Hilbert space'' }was denied by W. Johnson [4], who
constructed a series of examples of Banach spaces not isomorphic to the
Hilbert space, each of which has the hereditarily AP.

All these examples were ''sufficiently'' non symmetric and this fact allows
Johnson to pose a question (cf. [4]):

\textit{Whether there exists a Banach space }$X$\textit{\ with symmetric (or,
at least, subsymmetric) basis, distinct from the Hilbert space, with the
hereditarily approximation property?}

In this paper it will be obtained a partial answer on the Johnson's question.
Namely, it will be shown that there exists a Banach space $W$\ with a
subsymmetric basis, which is non equivalent to any symmetric basis, that has
the hereditarily AP.

Proofs are based on the notion of stable Banach spaces, due to J.-L. Krivine
and B. Maurey [5]. Namely, it will be shown that if a Banach space $X$ is not
isomorphic to a stable space then there exists a Banach space $W$ with a
spreading basis $\left(  w_{n}\right)  $, which is not equivalent to any
symmetric one, that is finitely representable in $X$.

Among examples presented in [4] there is a space $X_{J}$, defined as a
completition of vector space $c_{00}$ of all finite sequences of scalars (i.e.
sequences $\left(  a_{i}\right)  _{i=1}^{\infty}$, whose members all but a
finite number are vanished) by the norm
\[
\left\|  x\right\|  =\max\{\left\|  x\right\|  _{c_{0}},2^{-1}sup\left(
\sum\nolimits_{j=1}^{k_{n}}\left\|  A_{j}x\right\|  ^{2}\right)  ^{1/2}\},
\]
where $x=\left(  x_{n}\right)  _{n=1}^{\infty}=\left(  x_{n}\right)  $;
$A_{j}\subset\mathbb{N}$; $A_{j}\cap A_{k}=\varnothing$ ($j\neq k$); $A_{j}x $
denotes a sequence $\left(  y_{n}\right)  $ where $y_{n}=x_{n}$ for $n\in
A_{j}$; $y_{n}=0$ if $n\notin A_{j}$; the supremum is taken over all such
$A_{j}$'s that $\inf\{A_{j}:1\leq j\leq k_{n}\}\geq n$ and a sequence $\left(
k_{n}\right)  $ of natural numbers is chosen in a some special way.

As it was mentioned in [4], the space $X_{J}$ contains no subspaces that
isomorphic to the Hilbert space; it has the uniform approximation property and
enjoys the hereditarily AP. As it was noted before, since $X_{J}$ is not
isomorphic to any stable space, there exists a Banach space $W$ with a
spreading basis $\left(  w_{n}\right)  $ (since superreflexivity of $X_{J}$
the basis $\left(  w_{n}\right)  $ is subsymmetric), which is not equivalent
to any symmetric basis that is finitely representable in $X$. From the
author's paper [6] it$\operatorname*{rg}$ follows that $W$ has the
hereditarily AP. Moreover, $W$ has a more powerful property: any Banach space
that is finitely representable in $W$ has the AP.

The paper [6], mentioned above, was hardly compressed before its publication.
For this reason (and also for completeness), arguments from [6] will be
detailed here.

\section{Definitions and notations}

\begin{definition}
A Banach space $X$ is said to have the \textit{approximation property}%
\textbf{\ (}shortly\textbf{: }AP) if for every compact subset $K\subset X$ and
every $\varepsilon>0$ there exists a (bounded linear) operator $u:X\rightarrow
X$ with a finite dimensional range such that $\left\|  ux-x\right\|
\leq\varepsilon$ for all $x\in K$.

$X$ is said to have the $\lambda$-uniform approximation property ($\lambda
$-UAP) if there exists a function $f(n,\lambda):\mathbb{N\rightarrow N}$ such
that for any finite dimensional subspace $A$ of $X$ of dimension $\dim(A)=n$
there exists a finite rank operator $u:X\rightarrow X$, which is identical on
$A$, is of norm $\left\|  u\right\|  <\lambda$ and such that
$\operatorname*{rg}(u)=\dim(uX)\leq f(n,\lambda)$. If $X$ has the $\lambda
$-UAP for some $\lambda<\infty$, it will be said that $X$ has the UAP.
\end{definition}

\begin{definition}
It will be said that a Banach space $Y$ is finitely representable in a Banach
space $X$ (shortly $Y<_{f}X$ ) if for every $\varepsilon>0$ and every finite
dimensional subspace $A$ of $Y$ there exists a subspace $B$ of $X$ that is
$(1+\varepsilon)$-isomorphic to $A$ (i.e. there exists an isomorphism $u$
between $A$ and $B$ with $\mathcal{N}(u)=\left\|  u\right\|  \left\|
u^{-1}\right\|  <1+\varepsilon$).

$Y$ is crudely finite representable in $X$ ( in symbols: $Y<_{F}X$ ) if $Y$ is
isomorphic to a space $Z$ which is finitely representable in $X$
\end{definition}

\begin{definition}
Banach spaces $X$ and$\ Y$ are said to be finitely equivalent, shortly:
$X\sim_{f}Y$ (resp, crudely finitely equivalent, $X\sim_{F}Y$) if $X<_{f}Y$
and $Y<_{f}X$ (resp., if $X<_{F}Y$ and $Y<_{F}X$).

Any Banach space $X$ generates two classes:
\[
X^{f}=\{Y\in\mathcal{B}:X\sim_{f}Y\};\text{ \ }X^{F}=\{Y\in\mathcal{B}%
:X\sim_{F}Y\}.\text{\ }%
\]
\end{definition}

For any two Banach spaces $X$, $Y$ their \textit{Banach-Mazur distance }is
given by
\[
d(X,Y)=\inf\{\left\|  u\right\|  \left\|  u^{-1}\right\|  :u:X\rightarrow
Y\},
\]
where $u$ runs all isomorphisms between $X$ and $Y$ and is assumed, as usual,
that $\inf\varnothing=\infty$.

It is well known that $\log d(X,Y)$ defines a metric on each class of
isomorphic Banach spaces.

A set $\frak{M}_{n}$ of all $n$-dimensional Banach spaces, equipped with this
metric, is a compact metric space that is called \textit{the Minkowski
compact} $\frak{M}_{n}$.

The disjoint union $\cup\{\frak{M}_{n}:n<\infty\}=\frak{M}$ is a separable
metric space, which is called the \textit{Minkowski space}.

Consider a Banach space $X$. Let $H\left(  X\right)  $ be a set of all its
\textit{different} finite dimensional subspaces (\textit{isometric finite
dimensional subspaces of }$X$\textit{\ in }$H\left(  X\right)  $\textit{\ are
identified}). Thus, $H\left(  X\right)  $ may be regarded as a subset of
$\frak{M}$, equipped with the restriction of the metric topology of $\frak{M}$.

Of course, $H\left(  X\right)  $ need not to be a closed subset of $\frak{M}$.
Its closure in $\frak{M}$ will be denoted $\overline{H\left(  X\right)  }$.

From definitions it follows that $X<_{f}Y$ if and only if $\overline{H\left(
X\right)  }\subseteq\overline{H\left(  Y\right)  }$. Certainly, spaces $X$ and
$Y$ are \textit{finitely equivalent }( $X\sim_{f}Y$) if and only if
$\overline{H\left(  X\right)  }=\overline{H\left(  Y\right)  }$.

All spaces $Y$ from $X^{f}$ have the same set $\overline{H\left(  X\right)  }%
$. This set, uniquely determined by $X$ (or, equivalently, by $X^{f}$), will
be denoted by $\frak{M}(X^{f})$ and will be referred to as \textit{the
Minkowski's base of the class} $X^{f}$.

\begin{definition}
(Cf. $[7]$). For a Banach space $X$ its $l_{p}$-\textit{spectrum }$S(X)$ is
given by
\[
S(X)=\{p\in\lbrack0,\infty]:l_{p}<_{f}X\}.
\]
\end{definition}

Certainly, if $X\sim_{f}Y$ then $S(X)=S(Y)$. Thus, the $l_{p}$-spectrum $S(X)$
may be regarded as a property of the whole class $X^{f}$. So, notations like
$S(X^{f})$ are of obvious meaning.

\begin{definition}
Let $X$ be a Banach space. It is called \textit{superreflexive,} if every
space of the class $X^{f}$ is reflexive.
\end{definition}

Equivalently, $X$ is superreflexive if any $Y<_{f}X$ is reflexive.

\begin{definition}
Let $I$ be a set; $D$ be an ultrafilter over $I$; $\{X_{i}:i\in I\}$ be a
family of Banach spaces. An \textit{ultraproduct }$(X_{i})_{D}$ is a quotient
space
\[
(X)_{D}=l_{\infty}\left(  X_{i},I\right)  /N\left(  X_{i},D\right)  ,
\]
where $l_{\infty}\left(  X_{i},I\right)  $ is a Banach space of all families
$\frak{x}=\{x_{i}\in X_{i}:i\in I\}$, for which
\[
\left\|  \frak{x}\right\|  =\sup\{\left\|  x_{i}\right\|  _{X_{i}}:i\in
I\}<\infty;
\]
$N\left(  X_{i},D\right)  $ is a subspace of $l_{\infty}\left(  X_{i}%
,I\right)  $, which consists of such $\frak{x}$'s$\ $that
\[
\lim_{D}\left\|  x_{i}\right\|  _{X_{i}}=0.
\]
\end{definition}

If all $X_{i}$'s are all equal to a given $X\in\mathcal{B}$ then the
ultraproduct is said to be the \textit{ultrapower} and is denoted by $\left(
X\right)  _{D}$.

An operator $d_{X}:X\rightarrow\left(  X\right)  _{D}$ that sends any $x\in X$
to an element $\left(  x\right)  _{D}\in\left(  X\right)  _{D}$, which is
generated by a stationary family $\{x_{i}=x:i\in I\}$, is called the
\textit{canonical embedding }of $X$ into its ultrapower $\left(  X\right)
_{D}$.

It is well-known that a Banach space $X$ is finitely representable in a Banach
space $Y$ if and only if there exists such ultrafilter $D$ (over $I=\cup D$)
that $X$ is isometric to a subspace of the ultrapower $(Y)_{D}$.

Let $X$ be a Banach space.

\begin{definition}
A sequence $\{x_{n}:n<\infty\}$ of elements of $X$ is said to be
\end{definition}

\begin{itemize}
\item \textit{Spreading,} if for any $n<\infty$, any $\varepsilon>0$, any
scalars $\{a_{k}:k<n\}$ and any choosing of $i_{0}<i_{1}<...<i_{n-1}<...$;
$j_{0}<j_{1}<...<j_{n-1}<...$ of natural numbers
\[
\left\|  \sum\nolimits_{k<n}a_{k}x_{i_{k}}\right\|  =\left\|  \sum
\nolimits_{k<n}a_{k}x_{j_{k}}\right\|  .
\]

\item $C$-\textit{unconditional}, where $C<\infty$ is a constant, if
\[
C^{-1}\left\|  \sum\nolimits_{k<n}a_{k}\epsilon_{k}x_{i_{k}}\right\|
\leq\left\|  \sum\nolimits_{k<n}a_{k}x_{i_{k}}\right\|  \leq C\left\|
\sum\nolimits_{k<n}a_{k}\epsilon_{k}x_{i_{k}}\right\|
\]
for any choosing of $n<\infty$; $\{a_{k}:k<n\}$; $i_{0}<i_{1}<...<i_{n-1}<...$
and of signs $\{\epsilon_{k}\in\{+,-\}:k<n\}$.

\item \textit{Unconditional}, if it is $C$-unconditional for some $C<\infty$.

\item \textit{Symmetric,} if for any $n<\omega$, any finite subset
$I\subset\mathbb{N}$ of cardinality $n$, any rearrangement $\varsigma$ of
elements of $I$ and any scalars $\{a_{i}:i\in I\}$,
\[
\left\|  \sum\nolimits_{i\in I}a_{i}z_{i}\right\|  =\left\|  \sum
\nolimits_{i\in I}a_{\varsigma(i)}z_{i}\right\|  .
\]

\item \textit{Subsymmetric}, if it is both spreading and 1-unconditional.
\end{itemize}

Let $C<\infty$ be a constant. Two sequences $\{x_{n}:n<\infty\}$ and
$\{y_{m}:m<\infty\}$ are said to be $C$\textit{-equivalent} if for any finite
subset $I=\{i_{0}<i_{1}<...<i_{n-1}\}$ of $\mathbb{N}$ and for any choosing of
scalars $\{a_{k}:k<n\}$
\[
C^{-1}\left\|  \sum\nolimits_{k<n}a_{k}x_{i_{k}}\right\|  \leq\left\|
\sum\nolimits_{k<n}a_{k}y_{i_{k}}\right\|  \leq C\left\|  \sum\nolimits_{k<n}%
a_{k}x_{i_{k}}\right\|  .
\]

Two sequences $\left(  x_{n}\right)  $ and $\left(  y_{m}\right)  $ are said
to be \textit{equivalent} if they are $C$-equivalent for some $C<\infty$.

\section{Superstable classes of finite equivalence}

\begin{definition}
A Banach space $X$ is said to have the Tsirelson property if it does not
contain subspaces that are isomorphic to $l_{p}$ $(1\leq p<\infty)$ or $c_{0}$.
\end{definition}

The first example of a Banach space with such property was constructed by B.S.
Tsirelson [8].

From the other hand, there are classes $X^{f}$ that has ''anti-Tsirelson
property'': every representative of a such class contains some $l_{p}$. These
classes may be pick out by using the notion of \textit{stable Banach spaces},
which was introduced by J.-L. Krivine and B. Maurey [5].

\begin{definition}
A Banach space $X$ is said to be stable provided for any two sequences
$\left(  x_{n}\right)  $ and $\left(  y_{m}\right)  $ of its elements and
every pair of ultrafilter $D$, $E$ over $\mathbb{N}$
\[
\lim_{D\left(  n\right)  }\lim_{E\left(  m\right)  }\left\|  x_{n}%
+y_{m}\right\|  =\lim_{E\left(  m\right)  }\lim_{D\left(  n\right)  }\left\|
x_{n}+y_{m}\right\|  .
\]
\end{definition}

The notations $D(n)$ and $E(m)$ are used here (instead of $D$ and $E$) to
underline the variable ($n$ or $m$ respectively) in expressions like
$\lim_{D\left(  n\right)  }f(n,m)$.

\begin{definition}
(Cf. $[9]$). A Banach space $X$ is said to be superstable if every its
ultrapower $\left(  X\right)  _{D}$ is stable.
\end{definition}

\begin{theorem}
Let $X^{f}$ be a class of finite equivalence. $X^{f}$ contains a superstable
Banach space if and only if every space $Y\in X^{f}$ (and, as a consequence,
every space $W$, which is finitely representable in $X$) is stable.
\end{theorem}

\begin{proof}
If $Y\in X^{f}$ is superstable then each its subspace is stable because of the
property of a Banach space to be stable is inherited by its subspaces. Hence,
each subspace of every ultrapower $\left(  Y\right)  _{D}$ is stable as well,
because of $\{Z:Z<_{f}Y\}$ is coincide with the set
\[
\{Z:Z\text{ \ is isometric to a subspace of some ultrapower }\left(  Y\right)
_{D}\}.
\]

Conversely, if every $Y\in X^{f}$ is stable, then all ultrapowers of $Y$ are
stable too.
\end{proof}

\begin{definition}
A class $X^{f}$ of finite equivalence that contains a superstable space will
be called a superstable class.
\end{definition}

In [5] it was shown, that \textit{any stable Banach space }$X$\textit{\ is
weakly sequentially complete; every subspace of }$X$\textit{\ contains a
subspace isomorphic to some} $l_{p}$ ($1\leq p<\infty$).

\begin{definition}
(Cf. $[10]$) Let $X$ be a Banach space, $\left(  x_{n}\right)  \subset X$ be a
nontrivial normed sequence of elements of $X$ (i.e. $\left(  x_{n}\right)  $
contains no Cauchy subsequences); $D$ be an ultrafilter over $\mathbb{N}$. For
a finite sequence $\left(  a_{i}\right)  _{i=1}^{n}\subset\mathbb{R}^{n}$ let
\[
l\left(  \left(  a_{i}\right)  _{i=1}^{n}\right)  =\lim_{D\left(
m_{1}\right)  }\lim_{D\left(  m_{2}\right)  }...\lim_{D\left(  m_{n}\right)
}\{\left\|  \sum\nolimits_{k=1}^{n}a_{k}x_{m_{k}}\right\|  :m_{1}%
<m_{2}<...<m_{n}\}.
\]

Let $sm(X,\left(  x_{n}\right)  ,D)$ be the completition of the linear space
$c_{00}$ of all sequences $\left(  a_{i}\right)  _{i=1}^{\infty}$ of real
numbers such that all but finitely many $a_{i}$'s are equal to zero.

The space $sm(X,\left(  x_{n}\right)  ,D)$ is called a spreading model of the
space $X$, which is based on the sequence $\left(  x_{n}\right)  $ and on the
ultrafilter $D$.
\end{definition}

Clearly, any spreading model of a given Banach space $X$ has a spreading basis
and is finitely representable in $X$.

\begin{definition}
Let $X$ be a Banach space. Its $IS$-spectrum $IS(X)$ is a set of all
(separable) spaces $\left\langle Y,\left(  y_{i}\right)  \right\rangle $ with
a spreading basis $\left(  y_{i}\right)  $ which are finitely representable in
$X$.
\end{definition}

\begin{theorem}
A class $X^{f}$ is superstable if and only if every member $\left\langle
Y,\left(  y_{i}\right)  \right\rangle $ of its $IS$-spectrum has a symmetric basis.
\end{theorem}

\begin{proof}
Let $X^{f}$ be superstable; $Y\in X^{f}$. By [5], every spreading model of $Y
$ has a symmetric basis. Let $\left(  Y\right)  _{D}$ be an ultrapower by a
countably incomplete ultrafilter. Then (see [11]) $\left(  Y\right)  _{D}$
contains (as a subspace) every separable Banach space which is finitely
representable in $Y$. In particular, any space $\left\langle Z,\left(
z_{i}\right)  \right\rangle $ of $IS(X^{f})$ is isometric to a subspace of
$\left(  Y\right)  _{D}$.

Since $\left(  z_{i}\right)  $ is a spreading sequence, $sm(\left(  Y\right)
_{D},\left(  z_{i}\right)  ,E)$ is isometric to $\left\langle Z,\left(
z_{i}\right)  \right\rangle $ for any ultrafilter $E$. Hence, $\left(
z_{i}\right)  $ is a symmetric sequence.

Conversely, let every $\left\langle Z,\left(  z_{i}\right)  \right\rangle $
has a symmetric basis. Suppose that $X$ is not superstable. Then there exists
a space from $X^{f}$ which is not stable (it may be assumed that $X$ is not
stable itself). By [5] there are such sequences $\left(  x_{n}\right)  $ and
$\left(  y_{m}\right)  $ of elements of $X$ that
\[
\sup_{m<n}\left\|  x_{n}+y_{m}\right\|  >\inf_{m>n}\left\|  x_{n}%
+y_{m}\right\|  .
\]

Let $D$ be a countably incomplete ultrafilter over $\mathbb{N}$. Let ;
\[
X_{0}\overset{\operatorname{def}}{=}X\text{; \ }X_{n}\overset
{\operatorname{def}}{=}\left(  X_{n-1}\right)  _{D};\text{ }n=1,\text{
}2,\text{ }...;\text{ \ }X_{\infty}=\overline{\cup_{n\geq1}X_{n}}.
\]
Here is assumed that $X_{n}$ is a subspace of $X_{n+1}=\left(  X_{n}\right)
_{D}$ under the canonical embedding $d_{X_{n}}:X_{n}\rightarrow\left(
X_{n}\right)  _{D}$.

Let $D$, $E$ be ultrafilters over $\mathbb{N}$. Their product $D\times E$ is a
set of all subsets $A$ of $\mathbb{N}\times\mathbb{N}$ that are given by
\[
\{j\in\mathbb{N}:\{i\in\mathbb{N}:\left(  i,j\right)  \in A\}\in D\}\in E.
\]
Certainly, $D\times E$ is an ultrafilter and for every Banach space $Z$ the
ultrapower $\left(  Z\right)  _{D\times E}$ may be in a natural way identified
with $\left(  \left(  Z\right)  _{D}\right)  _{E}$.

So, the sequence $\left(  x_{n}\right)  \subset X$ defines elements
\begin{align*}
\frak{x}_{1}  &  =\left(  x_{n}\right)  _{D}\in\left(  X\right)  _{D};\\
\frak{x}_{2}  &  =\left(  x_{n}\right)  _{D\times D}\in\left(  \left(
X\right)  _{D}\right)  _{D};\\
&  ...\text{ \ \ \ \ }...\text{ \ \ \ \ }...\text{ \ \ \ \ }...\text{
\ \ \ \ }...\text{ \ \ \ \ }...\\
\frak{x}_{k}  &  =\left(  x_{n}\right)  \underset{k\text{ }times}%
{_{\underbrace{D\times D\times...\times D}}}\in\underset{k\text{ }%
times}{\underbrace{\left(  \left(  \left(  X\right)  _{D}\right)
_{D}...\right)  _{D}};}\\
&  ...\text{ \ \ \ \ }...\text{ \ \ \ \ }...\text{ \ \ \ \ }...\text{
\ \ \ \ }...\text{ \ \ \ \ }...
\end{align*}

Notice that $\frak{x}_{k}\in X_{k}\backslash X_{k-1}$. It is easy to verify
that $\left(  \frak{x}_{k}\right)  _{k<\infty}\subset X_{\infty}$ is a
spreading sequence. Since $X_{\infty}\in X^{f}$, it is symmetric. Moreover,
for any $z\in X$, where $X$ is regarded as a subspace of $X_{\infty}$ under
the direct limit of compositions
\[
d_{X_{n}}\circ d_{X_{n-1}}\circ...\circ d_{X_{0}}:X\rightarrow X_{n},
\]
the following equality is satisfied: for any pair $m$, $n\in\mathbb{N}$
\[
\left\|  x_{n}+z\right\|  =\left\|  x_{m}+z\right\|  .
\]
Since $\left(  x_{n}\right)  $ and $z$ are arbitrary elements of $X$, this
contradicts with the inequality $\sup_{m<n}\left\|  x_{n}+y_{m}\right\|
>\inf_{m>n}\left\|  x_{n}+y_{m}\right\|  $.
\end{proof}

This result may be generalized to classes of crudely finite equivalence.

\begin{definition}
A class $X^{F}$ of crudely finite equivalence is said to be crudely
superstable if it contains a superstable space.
\end{definition}

Certainly, any crudely superstable class $X^{F}$ has the property:

\textit{Every Banach space }$Y$\textit{, which is finitely representable in a
some space }$Z\in X^{F}$\textit{\ contains a subspace that is isomorphic to
some }$l_{p}$\textit{\ (}$1<p<\infty$\textit{)}.

\begin{theorem}
A class $X^{F}$ is crudely superstable if and only if for every $Z\in X^{F}$
its $IS$-spectrum $IS(Z)$ consists of spaces $\left\langle W,\left(
w_{n}\right)  \right\rangle $, whose natural bases $\left(  w_{n}\right)  $
are $c_{Z}$-equivalent to symmetric bases where the constant $c_{Z}$ depends
only on $Z$.
\end{theorem}

\begin{proof}
Let $X^{F}$ be crudely superstable. Then some $Y\in X^{F}$ is superstable and
for any $Z$, which is crudely finitely representable in $Y$, and for any space
$\left\langle W,\left(  w_{n}\right)  \right\rangle $ from the $IS$-spectrum
$IS(Z)$, its natural basis $\left(  w_{n}\right)  $ is $d(Z,Y_{1})$-equivalent
to symmetric one, where $Y_{1}\in Y^{f}$.

Conversely, let $X_{0}$ be such that every space $\left\langle W,\left(
w_{n}\right)  \right\rangle \in IS(X_{0})$ has a basis $\left(  w_{n}\right)
$, which is equivalent to a symmetric one (certainly, this is equivalent to
the assertion that for every $Z\in X^{F}$ its $IS$-spectrum $IS(Z)$ consists
of spaces $\left\langle W^{\prime},\left(  w_{n}^{\prime}\right)
\right\rangle $, whose natural bases $\left(  w_{n}^{\prime}\right)  $ are
equivalent to symmetric bases. It is easy to show that there exists a constant
$c_{X}$ such that every space $\left\langle W,\left(  w_{n}\right)
\right\rangle \in IS(X_{0})$ has a basis $\left(  w_{n}\right)  $ which is
$c_{X}$- equivalent to a symmetric one.

Indeed, let $\left(  c_{k}\right)  $ be a sequence of real numbers with a
property: for every $k<\infty$ there exists $\left\langle W_{k},\left(
w_{n}^{k}\right)  \right\rangle \in IS(X_{0})$ such that $\left(  w_{n}%
^{k}\right)  $ is $c_{k}$-equivalent to a symmetric basis and is not $c_{k-1}%
$-equivalent to any symmetric basis. Without loss of generality it may be
assumed that all spaces $\left(  W_{k}\right)  $ are subspaces of a space from
$\left(  X_{0}\right)  ^{f}$, e.g. $W_{k}\hookrightarrow X_{0}$. Consider an
ultrapower $\left(  X_{0}\right)  _{D}$ and its elements $\frak{w}_{k}=\left(
w_{n}^{k}\right)  _{D(k)}$. Clearly $\left(  \frak{w}_{k}\right)  _{k<\infty
}\subset\left(  X_{0}\right)  _{D}$ is a spreading sequence that is not
equivalent to any symmetric sequence.

Consider some $Z\in X^{F}$, such that $Z$ contains any space from its
$IS$-spectrum. Let $c_{Z}$ be the corresponding constant which was defined
above. Let $\{Z_{\alpha}:\alpha<\varkappa\}$ ($\varkappa$ is a cardinal
number) be a numeration of all subspaces of $Z$ that may be represented as
$\operatorname*{span}(z_{i}^{(\alpha)})$ (i.e. as a closure of linear span of
$\{:i<\infty\}$) for such sequences $\{z_{i}^{(\alpha)}:i<\infty\}$ (not
necessary spreading ones) that are $c_{Z}$-equivalent to symmetric sequences.

Using the standard procedure of renorming, due to A. Pe\l czy\'{n}ski [12], it
may be constructed a space $Z_{\infty}\approx Z$ such that $Z_{\infty}$
contains as a subspace every space $W_{\infty}$ from $IS(Z_{\infty})$, which
(by the renorming procedure) has a symmetric basis.

From the theorem 6 it follows that $Z_{\infty}$ (and, hence, the whole class
$\left(  Z_{\infty}\right)  ^{f}$) is superstable. Since $Z_{\infty}\in X^{F}$
the class $X^{F}$ is crudely superstable.

Indeed, it is sufficient to choose as a unit ball $B(Z_{\infty})=\{w\in
Z_{\infty}:\left\|  w\right\|  \leq1\}$ a convex hull of the union of a set
$B(Z)$ with sets $\{c_{Z}^{-1}j_{\alpha}B(i_{\alpha}Z_{\alpha}):$
$\alpha<\varkappa\}$, where $i_{\alpha}:Z_{\alpha}\rightarrow W_{\alpha}$ is
an isomorphism between $Z_{\alpha}$ and a space $W_{\alpha}$ with a symmetric
basis $\left(  w_{n}^{\left(  \alpha\right)  }\right)  $, which is given by
$i_{\alpha}\left(  z_{n}^{(\alpha)}\right)  =w_{n}^{\left(  \alpha\right)  }$
for all $n<\infty$; $\left\|  i_{\alpha}\right\|  \left\|  i_{\alpha}%
^{-1}\right\|  \leq c_{Z}$; $j_{\alpha}$ is an embedding of the unit ball
$B(W_{\alpha})$ in a set $c_{Z}B(Z_{\alpha})=\{z\in Z_{\alpha}:c_{Z}^{-1}z\in
B\left(  Z_{\alpha}\right)  \}$:
\[
B(Z_{\infty})=\operatorname*{conv}\{B(Z)\cup\left(  \cup\{c_{Z}^{-1}j_{\alpha
}B(i_{\alpha}Z_{\alpha}):\alpha<\varkappa\}\right)  \}.
\]
\end{proof}

\section{Factorization of operators}

Results of this section are based on some T. Figiel ideas [13].

Let $X,Y$ be Banach spaces, $L(X,Y)$ be the space of all (linear, bounded)
operators from $X$ to $Y$. Let $F(X,Y)$ be its subspace of all \textit{finite
rank operators}; $K(X,Y)$ be the space of all compact operators. Let $G(X,Y)$
be the closure of $F(X,Y)$ in the strong operator topology of $L(X,Y) $. $G$
will be called the space of \textit{approximated operators}.

\begin{definition}
Let $T\in L(X,Y)$. A pair $\left(  u,v\right)  $ of operators $u:X\rightarrow
Z $; $v:Z\rightarrow Y$ is called the factorization of $T$ through $Z$ if
$v\circ u=T$.
\end{definition}

A factorization $\left(  u,v\right)  $ of an operator $T$ is said to be a
$K$-factorization (resp., $G$-factorization) if both operators $u$ and $v$ are
compact (resp., approximated).

\begin{remark}
The approximation property may be expressed in the introduced terms.

Namely, a Banach space $X$ has the AP provided for every Banach space $Z$ the
following equality is fulfilled:
\[
K(Z,X)=G(Z,X).
\]
\end{remark}

\begin{remark}
According to $[13]$ if a Banach space $Z$ has the AP and for every Banach
space $X$ every operator $T:X\rightarrow Y$ has a $K$-factorization through
$Z$ then $Y$ also has the $AP$.

Indeed, let $T:X\rightarrow Y$; let $\left(  u,v\right)  $ be its
$K$-factorization through $Z$. Since $Z$ has the AP then $u\in G(Y,Z)$, i.e.
there is a sequence $\left(  u_{n}\right)  \subset F(X,Y)$ such that
$\lim\left\|  u-u_{n}\right\|  =0$. Certainly, $v\circ u_{n}\in G(X,Y)$;
\[
\left\|  T-v\circ u_{n}\right\|  =\left\|  v\circ u-v\circ u_{n}\right\|
\leq\left\|  v\right\|  \left\|  u-u_{n}\right\|  \rightarrow0.
\]

So, $T\in G(X,Y)$. Since $T$ and $X$ are arbitrary, result follows.
\end{remark}

\begin{definition}
A Banach space $X$ is said to be finitely decomposable if there exists a
constant $C_{X}$ such for any its subspace $X^{n}$ of finite codimension there
exists a space $Y$ such that $d(X,Y)<C$, which is finitely representable in
$X$. $C_{X}$ will be called the constant of finite decomposability of $X$.
\end{definition}

Immediately, if a class $X^{f}$ is generated by a finitely decomposable space
$X$ with the corresponding constant $C$ then every $Y\in X^{f}$ is also
finitely decomposable with the same constant $C_{Y}=C$. Such classes are also
called finitely decomposable classes.

\begin{theorem}
Let $X$ be a Banach space, $Y<_{f}X$ and $Y^{f}$ is finitely decomposable.
Then for every sequence $\frak{A}=\left(  A_{i}\right)  _{i<\infty}%
\subset\frak{M}\left(  Y^{f}\right)  $ every $Z\in X^{f}$ contains a subspace
$W_{\frak{A}}=W(\frak{A},Z)$ with the property:

\begin{itemize}
\item  for every $i<\infty$ there is a projection $P_{i}:W_{\frak{A}%
}\rightarrow W_{\frak{A}}$ such that
\[
\left\|  P_{i}\right\|  d(P_{i}W_{\frak{A}},A_{i})\leq C.
\]
\end{itemize}
\end{theorem}

\begin{proof}
Let $Z\in X^{f}$. Since $Y<_{f}X$ for every $\varepsilon>0$ and each
$i\in\mathbb{N}$ the space $A_{i}$ is $(1+\varepsilon)$-isomorphic to a
subspace of $Z$.

Let $i:A_{1}\rightarrow Z$ be a $(1+\varepsilon)$- isomorphic embedding; put
$iA_{1}=A_{1}^{\prime}$. There exist a finite codimensional subspace $Z^{1}$
of $Z$, which contains $A_{1}^{\prime}$, and a projection $P_{1}%
:Z^{1}\rightarrow A_{1}^{\prime}$,\ such that $\left\|  P_{1}\right\|
\leq1+\varepsilon$ (cf. [14]). Put $Q_{1}=Id_{Z^{1}}-P_{1}$. Put $Z^{\left(
1\right)  }=Q_{1}Z^{1}$.

Next we proceed by induction. Assume that $Z^{1}$, $Z^{2}$, ..., $Z^{n}$,
$Z^{\left(  1\right)  }$, $Z^{\left(  2\right)  }$, ..., $Z^{\left(  n\right)
}$ are already chosen. Since $Z^{\left(  n\right)  }$ is of finite codimension
in $Z$ and $Y^{f}$ is finite decomposable, $Z^{\left(  n\right)  }$ contains a
subspace $A_{n+1}^{\prime}$ that is $(1+\varepsilon^{n})$ isomorphic to
$A_{n+1}$. From [9] it follows that there exists a subspace $Z^{n+1}%
\hookrightarrow Z^{\left(  n\right)  }$ of finite codimension that contains
$A_{n+1}^{\prime}$ and a projection $P_{n+1}:Z^{n+1}\rightarrow A_{n+1}%
^{\prime}$ of norm $\left\|  P_{n+1}\right\|  \leq1+\varepsilon^{n}$. Put
$Q_{n+1}=Id_{Z^{n+1}}-P_{n+1}$ and $Z^{\left(  n+1\right)  }=Q_{n+1}Z^{n+1}$.

Clearly, the closure $W_{\frak{A}}$ of the linear span of the sequence
$\{A_{i}^{\prime}:i<\infty\}$ has the desired property.
\end{proof}

\begin{definition}
Let $Y$, $Z$ be Banach spaces. It will be said that $Z$ has the finite
projection property with respect to $Y^{f}$; shortly: $Z\in fpp(Y^{f})$
provided \textit{for a dense sequence }$\left(  A_{i}\right)  _{i<\infty
}\subset\frak{M}\left(  Y^{f}\right)  $\textit{\ and for every }$i<\infty
$\textit{\ there is a projection }$P_{i}:Z\rightarrow Z$\textit{\ such that
}$\left\|  P_{i}\right\|  d(P_{i}Z,A_{i})\leq C$.
\end{definition}

\begin{theorem}
Let $Y^{f}$ be a finite decomposable class; $Z\in fpp(Y^{f})$. Then for every
Banach space $X$ each approximated operator $T\in G(X,Y)$ admits
$G$-factorization through $Z$.
\end{theorem}

\begin{proof}
Let $T\in G(X,Y)$.

For every $i\in\mathbb{N}$ chose $S_{i}\in F(X,Y)$ such that $\left\|
T-S_{i}\right\|  \leq4^{-i}$. Let $S_{0}=0$; $T_{i}=S_{i}-S_{i-1}$. Since
\begin{align*}
\left\|  T_{i}\right\|   &  =\left\|  \left(  T-S_{i}\right)  -\left(
T-S_{i-1}\right)  \right\| \\
&  \leq\left\|  \left(  T-S_{i}\right)  \right\|  +\left\|  \left(
T-S_{i-1}\right)  \right\|  \leq2\cdot4^{-i+1},
\end{align*}
the series $\sum T_{i}$ absolutely converges to $T$.

Let $E_{i}=T_{i}X\hookrightarrow Y$. Since $Z\in fpp(Y^{f})$ and $Y^{f}$ is
finite decomposable, $Z$ contains a complemented subspace $W(\frak{A},Z)$ (by
the preceding theorem), which contains a complemented subspace $Y_{0}$,
isomorphic to the direct orthogonal sum $E_{0}=\sum\oplus E_{i}$.

Let $U:E_{0}\rightarrow Y_{0}$ be the corresponding isomorphism.

Put $E_{k}^{\prime}=UE_{k}\hookrightarrow Z$; $v_{k}=U\mid_{E_{k}}$.

Let $P:Z\rightarrow Y_{0}$ and $Q_{i}:Y_{0}\rightarrow E_{i}^{\prime}$ be
projections; $R_{i}=v_{i}^{-1}:E_{i}^{\prime}\rightarrow Y$ be isomorphic
embeddings. Put $R_{i}\circ Q_{i}\circ P=s_{i}$. Then $v_{i}\circ s_{i}%
=Q_{i}\circ P$ and $Q_{i}\circ P_{i}=Id_{E_{i}}$. Certainly,
\[
\left\|  s_{i}\right\|  \left\|  v_{i}\right\|  \left\|  v_{i}^{-1}\right\|
\leq c_{Y},
\]
where $c_{Y}$ depends only on $Y$.

Let $u\in G(X,Z)$ and $v\in G(Z,Y)$ be given by
\begin{align*}
ux  &  =\sum2^{i}v_{i}T_{i}x\text{ \ for }x\in X;\\
vz  &  =\sum2^{-i}s_{i}z\text{ \ \ \ for }z\in Z.
\end{align*}

Their norms may be estimated by
\begin{align*}
\left\|  u\right\|   &  \leq\sum2^{i}\left\|  v_{i}\right\|  \left\|
T_{i}\right\|  \leq\sum2^{i}\left\|  v_{i}\right\|  2\cdot4^{-i+1}\leq
4c_{Y};\\
\left\|  v\right\|   &  \leq\sum2^{-i}\left\|  s_{i}\right\|  \leq c_{Y}.
\end{align*}

So, $\left\|  u\right\|  \left\|  v\right\|  \leq4c_{Y}$.

The following equality shows that $\left(  u,v\right)  $ is a factorization of
$T$:
\[
v\circ u\left(  x\right)  =\sum2^{-i}s_{i}\sum2^{j}v_{j}T_{j}x=\sum s_{i}%
v_{i}T_{i}x=Tx.
\]

This proves the theorem.
\end{proof}

\begin{theorem}
Let $X$, $Y$, $Z$ be Banach spaces; $Z$ has the UAP; $Y$ is finitely
representable in $Z$ and the class $Z^{f}$ is finite decomposable. Then every
compact operator $T:X\rightarrow Y$ admits a $K$-factorization through a
subspace of $Z$.
\end{theorem}

\begin{proof}
Let $D$ be a such ultrafilter that $Y$ is isometric to a subspace of $\left(
Z\right)  _{D}$. Let $j:Y\rightarrow\left(  Z\right)  _{D}$ be the
corresponding isometry. Since $Z^{f}$ is finitely decomposable and $Y<_{f}Z$,
by the theorem 4 there exists a subspace $W$ of $\left(  Z\right)  _{D}$,
which has the finite projection property with respect to $Y^{f}$.

Let $T:X\rightarrow Y$ be a compact operator. Then $j\circ T\in K(X,\left(
Z\right)  _{D})$. Since $Z$ has the UAP, $\left(  Z\right)  _{D}$ also enjoys
this property (cf. [11]) and, hence, $j\circ T\in G(X,\left(  Z\right)  _{D}$.

By the theorem 5 there exists a $G$-factorization $\left(  u,v\right)  $ of
$j\circ T$ through $W$.

Let $W_{0}\hookrightarrow W$ be such that $uX\subseteq W_{0}\subseteq\left(
v^{-1}\circ j\right)  Y$. Put
\begin{align*}
\widetilde{u}\left(  x\right)   &  =u\left(  x\right)  \text{
\ \ \ \ \ \ \ \ \ \ \ \ \ for \ }x\in X;\\
\widetilde{v}\left(  z\right)   &  =\left(  j^{-1}\circ v\right)  \left(
w\right)  \text{ \ \ for \ }w\in W.
\end{align*}
Then $\left(  \widetilde{u},\widetilde{v}\right)  $ is a desired
$K$-factorization of $T$ through a subspace of $Z$.
\end{proof}

\section{A space with a subsymmetric basis that has the hereditarily AP}

\begin{theorem}
There exists a (superreflexive) Banach space $W$ with a subsymmetric basis
that is not equivalent to any symmetric one, such that every Banach space $Z$
which is crudely finite representable in $W$ has the approximation property.
\end{theorem}

\begin{proof}
Consider the Johnson's space $X_{J}$ (cf. [4]) that was described in the introduction.

Since $X_{J}$ has the Tsirelson property, it is not isomorphic to any stable
Banach space. By the theorem 3 there exists a space $\left\langle W,\left(
w_{n}\right)  \right\rangle \in IS(X_{J})$ which spreading basis $\left(
w_{n}\right)  $ is not equivalent to a symmetric one. Since $X_{J}$ is
superreflexive, $\left(  w_{n}\right)  $ is subsymmetric. Certainly, $W^{f}$
is finitely decomposable (with $c_{W}=1$).

Since $W<_{f}X_{J}$ it follows that there exists a subspace $W_{0}%
\hookrightarrow X_{J}$, which is finitely equivalent to $W$ and such that
$W_{0}\in fpp(W)$. Clearly, $W_{0}$ has the hereditarily AP and the uniform
approximation property (because of every subspace of $X_{J}$ enjoys UAP by
Johnson's construction [4]).

Let $Z<_{f}W\left(  \sim_{f}W_{0}\right)  $. Every compact operator
$U:Y\rightarrow Z$ for an arbitrary Banach space $Y$ admits a $K$%
-factorization $\left(  u,v\right)  $ through some subspace of $W_{0}$. Since
every subspace of $W_{0}$ has the AP, $Z$ also enjoys the AP.

Clearly, if $Z_{1}$ is isomorphic to $Z$ then it also has the AP.
\end{proof}

Since every subspace of $W$ is finitely representable in it, the desired
result is obtained.

\begin{theorem}
There exists a (superreflexive) Banach space $W$ with a subsymmetric basis
that is not equivalent to any symmetric one, such that every its subspace has
the approximation property.
\end{theorem}

\begin{problem}
Whether there exists a non-superreflexive Banach space that has the
hereditarily approximation property?
\end{problem}

Notice that if such a space would exist, it may be shown (in a similar way)
that there exists a nonreflexive space with a spreading basis that also enjoys
this property.

However, it is not known the answer on the following weaker question:

\begin{problem}
Whether there exists a non-superreflexive Banach space such that its $l_{p}%
$-spectrum consists of a single point $\{2\}$?
\end{problem}

\section{On the super quotient hereditarily approximation property}

The space $W$, constructed in the previous section, has a stronger
approximation property.

\begin{definition}
A Banach space $X$ has the quotient hereditarily approximation property
(shortly: QH-AP) provided each quotient of every its subspace has the
approximation property (AP).

$X$ has the super quotient hereditarily approximation property (super QH-AP)
if every Banach space $Y$, which is finitely representable in $X$ has the QH-AP.
\end{definition}

Recall that the aforementioned Johnson's space $X_{J}$ has the QH-AP and the UAP.

It will be needed the following result.

\begin{theorem}
Let $Y\in fpp(W)$ and $X$ be crudely finitely representable in $W$. Then there
exists an ultrapower $(Y)_{D}$, which contains a complemented subspace,
isomorphic to $X^{\ast\ast}.$
\end{theorem}

\begin{proof}
Consider the set $G(X^{\ast\ast})$ of all \textit{different} finite
dimensional subspaces of $X^{\ast\ast}$ (spaces $A$, $B\in G(X^{\ast\ast})$
are identified if and only if $A\cap B=A=B$), which will be indexed by
elements of a some set $J$: $G(X^{\ast\ast})=\{A_{j}:j\in J\}$.

Put $S_{j}=\{i\in J:A_{j}\hookrightarrow A_{i}\}$.

For any $n\in\mathbb{N}$ and any choosing $\{i_{1},i_{2},...,i_{n}\}\subset J$
there exists such $j\in J$ that
\[
S_{j}(\varepsilon)\subset S_{i_{1}}\cap S_{i_{2}}\cap...\cap S_{i_{n}}.
\]

Indeed, it is enough to chose $j\in J$ so that $A_{i_{k}}\hookrightarrow
A_{j}$ for all $k=1,2....,n$ .

Therefore a set $\{S_{j}:j\in J\}$ may be extended to a some ultrafilter $D$
over $J$.

Consider the ultrapower $(Y)_{D}$.

Obviously, there exists such $\lambda<\infty$ that for every finite
dimensional subspace $A$ of $X$ there are finite rank operators $T_{A}%
:A\rightarrow Y$ and $P_{A}:Y\rightarrow X$ such that $P_{A}\circ T_{A}%
\mid_{A}=Id_{A}$ and
\[
\left\|  T_{A}\right\|  \left\|  T_{A}^{-1}\right\|  \left\|  P_{A}\right\|
<\lambda.
\]
From the principle of local reflexivity it follows that the same is true for
subspaces of $X^{\ast\ast}$.

Let an embedding $V:X^{\ast\ast}\rightarrow(Y)_{D}$ be given by
\[
Vx^{\ast\ast}=(y_{i})_{D}\text{ \ for }x^{\ast\ast}\in X^{\ast\ast},
\]
where
\begin{align*}
y_{i} &  =T_{A_{i}}x^{\ast\ast}\text{ \ for \ }x^{\ast\ast}\in A_{i};\\
&  =0\text{ \ \ \ \ \ \ \ \ \ for \ }x^{\ast\ast}\notin A_{i}.
\end{align*}

Put $Q\left(  y_{i}\right)  _{D}=w^{\ast}-\lim_{D}P_{i}y_{i}^{\prime}$, where
$w^{\ast}-\lim_{D}$ means a limit by an ultrafilter $D$ in a weak* topology of
$X^{\ast\ast}$. Because of weakly* compactness of the unit ball of
$X^{\ast\ast}$ this limit exists.

Obviously, the operator $Q\circ V$ is identical on $X^{\ast\ast}$. The
operator $R=V\circ Q$ defines a projection from $\left(  Y\right)  _{D}$ onto
$X^{\ast\ast}$ of norm $\left\|  R\right\|  <\lambda$.
\end{proof}

\begin{theorem}
There exists a superreflexive Banach space non isomorphic to the Hilbert space
which has a subsymmetric basis and enjoys the super QH-AP.
\end{theorem}

\begin{proof}
Let $W$ be the space with a subsymmetric basis, which has the hereditarily AP
and the UAP that was constructed in the theorem 7. Let $X_{J}$ be the
Johnson's space that has the QH-AP. Let $W_{0}$ be a subspace of $X_{J}$, such
that $W_{0}\sim_{f}W$ and $W_{0}\in fpp(W)$. It exists by the theorem 4. Since
$W_{0}\hookrightarrow X_{J}$, it enjoys the QH-AP and the UAP as well.

Consider its conjugate $(W_{0})^{\ast}$. It has the hereditarily AP as well as
UAP and, by the theorem 7, enjoys the super hereditarily AP (i.e., every space
that is finitely representable in $(W_{0})^{\ast}$ enjoys the hereditarily AP.

Let $Z$ be crudely finitely representable in $W$ (and, hence, in $W_{0}$). Let
$Y$ be a quotient of a subspace $V$ of $Z$. Surely, $Y^{\ast}$ may be
identified with a subspace of $V^{\ast}$. Since $W_{0}\in fpp(W)$, by the
preceding theorem there exists such ultrapower $(W_{0})_{D}$, which contains a
complemented subspace, isomorphic to a given space, which is crudely finitely
representable in $W.$ So, $V$ may be regarded as an isomorph of a complemented
subspace of $(W_{0})_{D}.$Hence, $V^{\ast}$ is crudely finite representable in
$(W_{0})^{\ast}$, and, as a consequence, has the hereditarily AP.This yields
that $Y^{\ast}$, as a subspace of $V^{\ast}$, has the AP.

Since $Z$ is arbitrary, this proves the theorem.
\end{proof}

\section{References}

\begin{enumerate}
\item  Enflo P. \textit{A counterexample to the approximation problem in
Banach spaces}, Acta Math. \textbf{130} (1973) 309-317

\item  Szankowski A. \textit{Subspaces without approximation property}, Israel
J. Math. \textbf{30} (1978) 123-129

\item  Reinov O.I. \textit{Banach spaces without the approximation property
}(in Russian)\textit{,} Funct. Analysis and Appl. \textbf{16:4} (1982) 84-85

\item  Johnson W.B. \textit{Banach spaces all of whose subspaces have the
approximation property}, Special topics in Appl. Math. (1980) 15-26

\item  Krivine J.-L., Maurey B. \textit{Espaces de Banach stables}, Israel J.
Math. \textbf{39} (1981) 273-295

\item  Tokarev E.V. \textit{Banach spaces with the superapproximation
property} (transl. from Russian), Ukrainian Mathematical Journal \textbf{38}
(1986) 230-231

\item  Schwartz L. \textit{Geometry and probability in Banach spaces}, Bull.
AMS \textbf{4:2} (1981) 135-141

\item  Tsirelson B.S. \textit{Not every Banach space contains }$l_{p}%
$\textit{\ or }$c_{0}$ (in Russian)\textit{,} Funct. Analysis and Appl.
\textbf{8:2} (1974) 57-60

\item  Raynaud Y. \textit{Espaces de Banach superstables}, C. R. Acad. Sci.
Paris, S\'{e}r. A. \textbf{292:14} (1981) 671-673

\item  Brunel A., Sucheston L. \textit{On }$\mathit{B}$\textit{-convex Banach
spaces}, Math. System Theory \textbf{7} (1973) 294-299

\item  Heinrich S. \textit{Ultraproducts in Banach space theory}, J. Reine
Angew. Math. \textbf{313} (1980) 72-104

\item  Pe\l czy\'{n}ski A. \textit{Projections in certain Banach spaces},
Studia Math. \textbf{19} (1960) 209 - 228

\item  Figiel T. \textit{Factorization of compact operators and applications
to the approximation problem}, Studia Math. \textbf{45 }(1973) 191-210

\item  Gurarii V.I. On spreadings and inclinations of subspaces of Banach
spaces Theory of Funct., Funct. Analysis and Appl.\textbf{\ 1 }(1965) 194-204
\end{enumerate}
\end{document}